\documentclass[11pt]{article}

\usepackage{amsmath,amsfonts}
\usepackage{color}
\definecolor{Blue}{rgb}{0.1,0.1,0.9}
\definecolor{Red}{rgb}{0.9,0.1,0.1}

\newcommand{\ex}{\mathbf{E}}
\newcommand{\pr}{\mathbf{P}}
\newcommand{\R}{\mathbb{R}}

\newcommand{\EE}{\mathcal{E}}
\newcommand{\Rd}{{\R^d}}
\newcommand{\Rt}{{\R^2}}
\newcommand{\M}{{\mathcal{M}}}

\newcommand{\N}{{\mathcal{N}}}
\newcommand{\F}{{\mathcal{F}}}
\newcommand{\LL}{{L}}
\newcommand{\CC}{{C_c}}
\newcommand{\Lp}{{\LL^p}}
\newcommand{\Lq}{{\LL^q}}
\newcommand{\Li}{{\LL^\infty}}
\newcommand{\Lo}{{\LL^1}}
\newcommand{\Lr}{{\LL^r}}
\newcommand{\Lt}{{\LL^2}}

\newcommand{\bC}{\mathbf{C}}

\newcommand{\ind}{\mathbf{1}}
\newcommand{\wt}{\widetilde}
\newcommand{\pl}{\|}

\newcommand{\proof}{\noindent \textsc{Proof: }}

\newcommand{\qed}{\ensuremath{\hspace{10pt}\square}}

\newtheorem{theorem}{Theorem}

\newtheorem{lemma}{Lemma}

\newtheorem{example}{Example}

\DeclareMathOperator{\supp}{supp}

\begin{document}
\sloppy

%                            ---------- o ----------

\begin{titlepage}
\title{\bf L\'evy processes and Fourier multipliers}
\author{Rodrigo Ba\~nuelos\thanks{Supported in part by NSF Grant
\#  0603701-DMS}\\Department of  Mathematics
\\Purdue University\\West Lafayette, IN
47906\\banuelos@math.purdue.edu\and  Krzysztof Bogdan
\thanks{Supported in part by grant DBN 1 P03A 026 29 }
\\ Institute of Mathematics\\
 Wroc{\l}aw University of Technology\\
50-370 Wroc{\l}aw, Poland\\bogdan@im.pwr.wroc.pl
} \maketitle

%%                            ---------- o ----------
\footnotetext{2000 {\it MS Classification}\/:
42B15, 60G51 (Primary), 60G46, 42B20 (Secondary).}

\begin{abstract}
We study Fourier multipliers which result from modulating jumps of
L\'evy processes. Using the theory of martingale transforms we prove
that these operators are bounded in $L^p(\Rd)$ for $1<p<\infty$ and we
obtain the same explicit bound for their norm as the one known 
for the second order Riesz transforms.   
\end{abstract}
\end{titlepage}

%\footnotetext{2000 {\it MS Classification}\/:
%42B15, 60G51 (Primary), 60G46, 42B20 (Secondary).\\
%{\it Key words and phrases}\/: Fourier multiplier, process
%with independent increments, 
%martingale transform, singular integral}
%Legend: Within 42Bxx Fourier analysis in several variables we have:
%42B15  Multipliers
%42B20  Singular integrals (Calder?n-Zygmund, etc.)
%60G51  Processes with independent increments
%60G46  Martingales and classical analysis
%60G44  Martingales with continuous parameter

%                            ---------- o ----------

\section{Introduction}
One of the most basic examples of 
 Calder\'on--Zygmund singular
integrals in $\Rd$ is the collection of Riesz transforms
(\cite{stein}), 
%on $\Rd$ are defined by for any $f\in L^p(\Rd)$ by 
\begin{displaymath}
R_jf(x)=\frac{\Gamma\left(\frac{d+1}{2}\right)}{\pi^{(d+1)/2}}
\;p.v.\int_{\Rd} {x_j-y_j\over |x-y|^{d+1}}f(y)\, dy\,,\quad
j=1, 2, \dots, d\,.   
\end{displaymath}
They can be represented as Fourier
multipliers with symbols $i\xi_j/|\xi|$,
%of the type considered by Marcinkiewicz: 
\begin{displaymath}
\widehat{R_jf}(\xi)=\frac{i\xi_j}{|\xi|}\hat f(\xi)\,,\quad  f\in L^2(\Rd)\,.
%  f\in L^2(\Rd),  j=1, 2, \dots, d.  
\end{displaymath}
%Iterating this one more time gives 
Therefore the second order Riesz transforms $R_j^2$ satisfy
% the multiplier operators 
\begin{equation}\label{eq:rt2}
\widehat{R_j^2f}(\xi)=-\frac{\xi_j^2}{|\xi|^2}\hat f(\xi)\,,\quad
%  f\in L^2(\Rd), 
j=1, 2, \dots, d\,.  
\end{equation}
%and 
%\begin{equation}
%\widehat{R_j R_kf}(\xi)=-\frac{\xi_j\xi_k}{|\xi|^2}\hat f(\xi),  f\in L^2(\Rd), j, k=1, 2, \dots, d.  
%\end{equation}
%These operators are Calder\'on--Zygmund singular integrals.  
It follows from the general theory of singular integrals (see Stein \cite{stein}) that there
exist constants $C_p$ and $C_p'$ such that $\|R_jf\|_p\leq
C_p\|f\|_p$ and $\|R_j^2 f\|_p\leq C_p'\|f\|_p$ for every
$1<p<\infty$.  There has been considerable interest in recent years in
obtaining the best values for these constants.  It
was shown in  \cite{iwaniec-martin} that 
 \begin{equation}\label{cot}
 \|R_jf\|_p \leq \cot \left({\pi\over 2p^*}\right)\|f\|_p\,,\quad f\in L^p(\Rd)\,, 
\end{equation} 
and that $\cot\left({\pi\over 2p^*}\right)$ is the best (smallest) possible
constant for this inequality.  Here and below,
 \begin{equation}
  \label{eq:defp}
1<p<\infty\,,\quad q=p/(p-1)\,,\quad p^*=\max(p, q)\,, 
\end{equation}
so that 
$$p^*-1=\max\{p-1, (p-1)^{-1}\}.$$
An alternative proof of (\ref{cot}) is given in \cite{banuelos-wang}
by applying the martingale transform techniques of Burkholder
(\cite{Bur8}, \cite{burkholder}) to stochastic integrals obtained from
composing harmonic functions with  Brownian motion.    
Using a similar approach,  it is
also proved in \cite{banuelos-wang} that 
\begin{equation}\label{riesz1}
\|R_jR_kf\|_p\leq (p^*-1)\|f\|_p,  \quad j\not = k\,
\end{equation}
and that   
\begin{equation}
  \label{eq:nrt}
\|R_j^2f\|_p\leq (p^*-1)\|f\|_p\,.
% for all $j=1, \dots, d$. 
\end{equation}
%By applying the techniques of \cite{banuelos-wang} 
%to space time martingales obtained by composing 
%caloric functions (solutions of the heat equation)
% with space time Brownian motion, it is shown in
%\cite{banuelos-mendez} that 
%%if $f\in L^p(\Rd)$, $1<p<\infty$, then
%\begin{equation}
%  \label{eq:c2rt}
%\|\sum_{j=1}^d a_jR_j^2f\|_p\leq (p^*-1)\|f\|_p\,,
%\end{equation}
%if $a_j\in \{-1, 0, 1 \}$, and
%\begin{equation}
%  \label{eq:i2rt}
%\|2R_jR_kf\|_p\leq (p^*-1)\|f\|_p\,, \quad j\not = k\,.
%\end{equation}
The above operators are closely related to the Beurling--Ahlfors operator 
\begin{equation}
  \label{eq:ba}
  B f(z)=-\frac{1}{\pi}\, p.v. \int_\mathcal{C}
  \frac{f(\zeta)}{(z-\zeta)^2}d\zeta_1 d\zeta_2\,.
\end{equation} 
Indeed, the Beurling--Ahlfors operator is a singular integral of
even kernel whose Fourier multiplier is $\overline{\xi}/\xi$ and hence
$B=-R^2_1+R^2_2+2i R_2 R_1$ (see \cite{banuelos-mendez}). 
The computation of the norm of $B$ on $L^p(\bC)$
has been a problem of considerable interest for many years now. 
In  \cite{Leh1}, Lehto showed that  $\|B\|_p \geq  p^*-1$ and T. Iwaniec
conjectured in \cite{iwaniec} that $\|B\|_p = p^*-1$.
%\begin{equation}\label{iwaniec-conj} 
%\|B\|_p = p^*-1,
%\end{equation}

%Up to now, all the techniques which give explicit estimates on the norm
%of $B$ depend heavily on the work of Burkholder who succeeds in
%computing best constants of martingale transforms (\cite{Bu1},
%\cite{Bu5}). 
In \cite{banuelos-wang}, the 
martingale inequalities of Burkholder, together with the representation of the $B$ as a conditional expectation of certain stochastic integrals, were used to prove the bound $\|Bf\|_p\leq
4(p^*-1)\|f\|_p$, for general complex valued $f$ and that  $\|Bf\|_p\leq
2\sqrt{2}(p^*-1)\|f\|_p$, for for real valued functions $f$.  In \cite{NV}
Nazarov and Volberg improved the bound to $2(p^*-1)$ for general $f$ and to $\sqrt{2}(p^*-1)$ for real valued $f$
using an 
analytic (Littlewood--Paley inequalities) approach with Bellman functions that also rests on 
the martingale inequalities of Burkholder. A different proof of the Nazarov-Volberg bounds 
was given in
\cite{banuelos-mendez} using essentially the same proof as the one in
\cite{banuelos-wang} but applied to space time Brownian martingales.
In \cite{volberg-oliver}, Dragi\v cevi\'c and Volberg  refined the
Nazarov-Volberg  techniques and obtained that for general $f$, 
\begin{equation}\label{vol-dra1} 
\|Bf\|_p\leq \sqrt{2}(p-1){\left(\frac{1}{2\pi}\int_0^{2\pi}|\cos(\theta)|^p
d\theta\right)}^{-\frac{1}{p}}\|f\|_p, \hspace{5mm} 2\leq p<\infty, 
\end{equation}
and that for real valued $f$, 
\begin{equation}\label{vol-dra2} 
\|Bf\|_p\leq (p-1){\left(\frac{1}{2\pi}\int_0^{2\pi}|\cos(\theta)|^p
d\theta\right)}^{-\frac{1}{p}}\|f\|_p, \hspace{5mm} 2\leq p<\infty.
\end{equation}

By a further refinement of the techniques in  \cite{banuelos-mendez}, 
it is proved in \cite{banuelos-janakiraman} that 
\begin{equation}\label{main}
\|Bf\|_p\leq \sqrt{2(p^2-p)}\,{\|f\|}_p, \hspace{5mm} 2\leq p<\infty,
\end{equation}
for general complex valued $f$, and that 

\begin{equation}\label{real_result}
\|Bf\|_p\leq \sqrt{p^2-p}\,{\|f\|}_p,\hspace{5mm} 2\leq p<\infty,
\end{equation}
for real valued $f$.  Dividing both bounds in (\ref{vol-dra1}) and (\ref{main}) by $p$ and letting this go to infinity both give $\sqrt{2}$.  However, asymptotically the estimate (\ref{main}) is slightly better as can be easily checked.  Interpolation and the bound in (\ref{main}) gives the general bound 
$\|B\|_p\leq 1.575(p^*-1)$ for the norm of the operator.  
For more information on Iwaniec's conjecture and its connections to
quasiconformal mappings and other areas of nonlinear PDE, we refer the
reader to 
\cite{iwaniec-martin}, \cite{iwaniec}, \cite{banuelos-wang},
\cite{NV},  \cite{volberg-oliver}, \cite{banuelos-mendez},
\cite{banuelos-janakiraman}. 

 The purpose of the present paper is to explore martingale techniques
to study Fourier multipliers which arise when the Brownian motion used
to define stochastic integrals leading to the Riesz
transforms
is replaced by the more general symmetric
L\'evy process. This leads to a large family of multipliers which
generalize the second order Riesz transforms. We obtain the upper bound
$p^*-1$ for their norms in $L^p(\Rd)$, which is the best known to date
in the case of the second order Riesz transforms $R_j^2$.

Let $ V\geq 0$ be a L\'evy measure on $\Rd$, that is $ V(\{0\})=0$,
$ V\neq 0$, and
\begin{equation}
  \label{eq:clm}
\int_\Rd \min(|x|^2,1) V(dx)<\infty\,.
\end{equation}
Assume that $ V$ is symmetric: $ V(-B)= V(B)$. Let
$\phi$ be complex-valued, Borel measurable and symmetric: $\phi(-z)=\phi(z)$, and
assume that
$$
|\phi(z)|\leq 1\,,\quad z\in \Rd\,.
$$
%Let
%\begin{equation}
%  \label{eq:defp}
%1<p<\infty\,,\quad q=p/(p-1)\,,\quad p^*=\max(p, q)\,, 
%\end{equation}
%so that $p^*-1=\max\{p-1,(p-1)^{-1}\}$.
\begin{theorem} \label{th:bfm}
  The Fourier multiplier with the symbol 
  \begin{equation}
    \label{eq:fs}
    M(\xi)=\frac{\int_\Rd(\cos\xi\!\cdot\! z-1)\phi(z) V(dz)}{\int_\Rd (\cos
    \xi\!\cdot\! z-1) V(dz)}\,,
  \end{equation}
is bounded on $\Lp(\Rd)$ for $1<p<\infty$, with the norm at most $p^*-1$.
% I think we should really say here what we mean
%by this.  It adds a little more information. 
That is, if we define the operator $\M$ on $L^2(\Rd)$ by 
$$
\widehat{\M f} (\xi)=M(\xi)\hat f(\xi),
$$
then $\M$ has a unique bounded linear extension to $L^p(\Rd)$, $1<p<\infty$, and 
\begin{equation}
  \label{eq:dM}
\|\M f\|_p\leq (p^*-1)\|f\|_p.
\end{equation}
\end{theorem}

We note that 
the boundedness of our multipliers on $L^p(\Rd)$ does not follow directly
from the H\"ormander multiplier theorem (\cite{stein}, page 96) because their
symbols (\ref{eq:fs}) generally lack sufficient differentiability.
However, for certain special cases (such as those mentioned in
(\ref{eq:pm}) below), general $L^p$ bounds can be obtained from the
Marcinkiewicz multiplier theorem, see Stein \cite{stein}, page 109.

%{\color{red}
As we already mentioned, the technique used here consists in representing
singular integrals and other Fourier multiplies by means of conditional
expectations of stochastic integrals. This approach has origins in the
paper of Gundy and Varopoulos \cite{gunvor},
and was widely applied to stochastic integrals based on
standard Brownian motion or the space-time Brownian motion
(see, e.g., \cite{banuelos-mendez} and \cite{banuelos-janakiraman}). 
As it is well-known, the Brownian motion at the
times when its last coordinate first reaches a certain level is a
space-time Cauchy process. 
McConnell studied in \cite{mcc} 
the resulting Cauchy process 
and related martingales (called parabolic martingales below) by a
discretization method (\cite[(3.9)]{mcc}),
to extend the classical H\"ormander multiplier theorem to functions 
taking values in Banach spaces with the unconditional martingale
difference (UMD) sequence property. 
The study can be considered a precursor of our development (see also
\cite{banuelos-mendez}).
%To the best of our knowledge, the present paper is the first application
%of more general L\'evy processes to obtain similar results.
In the present paper we employ an integral representation of parabolic martingales, 
and results of \cite{MR1334160} to obtain explicit estimates for the
$L^p$ norms of the considered multipliers, our main goal for this
study. 
%}

A word about our notation. We always assume Borel measurability of
considered sets and functions below.
By $\Lr=L^r(\Rd)$, with $1\leq r<\infty$,  we will denote the set of
complex-valued functions $g$ such that
$$
\pl g\pl_r =\left[\int_\Rd |g(x)|^r dx\right]^{1/r}<\infty\,,
$$
$\Li$ are those $g$ for which $\pl g\pl_\infty =\sup_{x\in \Rd}
|g(x)|<\infty$, 
and $\CC$ consists of continuous compactly supported functions $g$.
For $g\in \Lo$ its Fourier transform is defined as
\begin{displaymath}
  \F g(\xi)=\widehat{g}(\xi)=\int_\Rd e^{i \xi\cdot z}g(z)dz\,,\quad \xi\in \Rd\,. 
\end{displaymath}
By Plancherel's Theorem $\pl\widehat{g}\pl_2=(2\pi)^d\pl g\pl_2$, and
$\F$ extends to a continuous linear bijection of $\Lt$. 
Thus the Fourier multiplier $\M$ of Theorem~\ref{th:bfm} has the norm on $\Lt$ equal to
$\|M\|_\infty\leq 1$, see (\ref{eq:fs}).
Theorem~\ref{th:bfm} will be proved by verifying that (\ref{eq:dM}) holds
for every $f\in \CC$. This yields that $\M$ has a unique bounded
linear extension to $\Lp$, denoted also $\M$, satisfying (\ref{eq:dM}) for every $f\in \Lp$. 
% $1<p<\infty$. 

To give an example, let $\alpha\in (0,2)$ and $j=1,\ldots,d$. 
We have that (\ref{eq:dM}) holds when the multiplier has the symbol
\begin{equation}
  \label{eq:pm}
  M(\xi)=\frac{|\xi_j|^\alpha}{|\xi_1|^\alpha+\cdots+|\xi_d|^\alpha}\,,\quad
  \xi=(\xi_1,\ldots,\xi_d)\in \Rd\,.
\end{equation}
%H\"ormander's
%theorem (\cite{stein}) cannot be used to prove even the mere boundedness of such
%multiplier in $\Lp$ for $p\neq 2$ for general $\alpha$ and $d$.
We also note that (\ref{eq:dM}) extends to multipliers
whose symbols may be obtained as pointwise limits of symbols of the
from (\ref{eq:fs}). For instance, (\ref{eq:nrt}) can be obtained  by
letting $\alpha\to 2$ in (\ref{eq:pm}), see (\ref{eq:rt2}).

Here is the composition of the paper.
The proof of Theorem~\ref{th:bfm} is given in Section~\ref{s:p}.
In Section~\ref{s:m} we make some additional consideration, for example we
examine (\ref{eq:pm}).
The paper is essentially self-contained except for the $\Lp$
estimates for  differentially subordinate martingales, which in our case follow form 
the work of G. Wang  \cite{MR1334160}. 

%                            ---------- o ----------
%\section{L\'evy system}
\section{Proof of Theorem~\ref{th:bfm}}\label{s:p}
We first describe the setup which will be used in the proof of
the result. Let $\nu\geq 0$ be a {\it finite}\/ measure on $\Rd$ not
charging the origin. Assume that $\nu$ is symmetric: $\nu(-B)=\nu(B)$,
and $|\nu|=\nu(\Rd)>0$. Let $\wt{\nu}=\nu/|\nu|$. 
Let $\pr$ and $\ex$ be the probability and expectation for a family of
independent random variables $T_i$ and $Z_i$, $i=\pm 1,\pm 2,\ldots$, where each $T_i$ is
exponentially distributed with $\ex T_i=1/|\nu|$, and each $Z_i$ has
$\wt{\nu}$ as the distribution.
We let $S_i=T_1+\cdots+T_i$ for $i=1,2,\ldots$, and
$S_i=-(T_{-1}+\cdots+T_{i})$ for $i=-1,-2,\ldots$. 
For $-\infty<s<t<\infty$ we let $X_{s,t}=\sum_{s< S_i\leq t} Z_i$, and
$X_{s,t-}=\sum_{s< S_i< t} Z_i$.
We note that $\N(B)=\#\{i:\;(S_i,Z_i)\in B\}$ is a Poisson random measure on
$\R\times \Rd$ with intensity measure $dv\,\nu(dx)$, and
$X_{s,t}=\int_{s< v\leq t} x \N(dvdx)$ is the L\'evy-It\^o
decomposition of $X$ (\cite{MR1739520}).
Let $N(s,t)=\N((s,t]\times \Rd)$ be the number of {\it signals} $S_i$
such that $s<S_i\leq t$. 

For the reader's convenience we give an elementary proof of what
amounts to the L\'evy system for $X$ (see \cite[VII.68]{MR745449} for
more general results).
\begin{lemma} \label{lem:ls}
If 
{the Borel measurable} function  $F\,:\; \R\times \Rd\times \Rd\to\R$ is either 
  nonnegative or bounded, and $s\leq t$, then  
  \begin{equation}
    \label{eq:els}
    \ex \sum_{s< S_i\leq t} F(S_i,X_{s,S_{i}-}, X_{s,S_i})
= \ex \int_s^t \int_\Rd F(v,X_{s,v-},
X_{s,v{\color{blue}-}}+z)\nu(dz)dv\,.
  \end{equation}
\end{lemma}
\proof
Since the arrival time of the $n$-th signal has the gamma
distribution, 
%we have
\begin{eqnarray*}
  LHS &=& \sum_{-\infty<i<\infty} \ex\left\{ F(S_i,X_{s,S_{i}-}, X_{s,S_i})
  \ind_{s< S_i \leq t}\right\}\\
  &=& \sum_{n=0}^\infty \int_\Rd \int_\Rd \int_s^t
  F(v,y,y+z)\frac{|\nu|^{n+1}(v-s)^n}{n!}e^{-|\nu|(v-s)}dv\wt{\nu}^{*n}(dy)\wt{\nu}(dz)\\
  &=& \int_s^t\int_\Rd \int_\Rd 
  F(v,y,y+z)e^{-|\nu|(v-s)}e^{*(v-s)\nu}(dy)\nu(dz)dv\,.  
 \end{eqnarray*}
Here $\mu^{*n}$ is the $n$-fold convolution of a measure $\mu$ and 
$e^{*\mu}=\sum_{n=0}^\infty \mu^{*n}/n!$ denotes the convolution exponent
of $\mu$. 
In what follows we will use the following two well-known facts.
\begin{enumerate}
\item \noindent
First, conditionally
on $N(s,t)=n$, the consecutive signals in $(s,t]$
are uniformly distributed on $\{(s_1,\ldots,s_n)\,:\;s< s_1
\leq \ldots \leq s_n\leq t\}$.

\noindent
\item Second, let $s< v \leq t$. Let $T g(v)=\int_v^t g(u)du$ for measurable and bounded or
nonnegative function $g$. By induction, for $n=1,2,3,\ldots$,
$$T^{n}g(v)=T(T^{n-1} g)(v)=\frac{1}{(n-1)!}\int_v^t g(u)(u-v)^{n-1}du\,.$$
%and, similarly, if $S g(v)=\int_s^v g(u)du$, then  
%$$S^{n} g(v)=\frac{1}{(n-1)!}\int_s^v g(u)(v-u)^{n-1}du\,.$$
\end{enumerate}
We have
\begin{eqnarray*}
  RHS
&=&  \sum_{n=0}^\infty \ex \left\{\int_s^t \int_\Rd 
  F(v,X_{s,v-}, X_{s,v-}+z)\nu(dz)dv|N(s,t)=n\right\}
  \frac{|\nu|^{n}(t-s)^n}{n!}e^{-|\nu|(t-s)}\\
&=&  \sum_{n=0}^\infty
\frac{|\nu|^{n}(t-s)^n}{n!}e^{-|\nu|(t-s)}\frac{n!}{(t-s)^n}
\int_\Rd \int_\Rd \int_s^t ds_1 \int_{s_1}^t ds_2 \ldots \int_{s_{n-1}}^t ds_n  \\
&&\sum_{k=0}^n \int_{s_k}^{s_{k+1}}
 F(v,y,y+z)dv\wt{\nu}^{*k}(dy)\nu(dz)\,,
\end{eqnarray*}
where $s_0=s$ and $s_{k+1}=t$ for $k=n$.
Changing notation involving $v$ and $s_k$ we obtain
\begin{eqnarray*}
&&\!\!\!\!\!\!\!\!\!\!\!\!\!\!\!  \sum_{n=0}^\infty
|\nu|^{n}e^{-|\nu|(t-s)}
\int_\Rd \int_\Rd \sum_{k=0}^n 
\int_s^t ds_1 \ldots \int_{s_k}^t ds_{k+1}
 F(s_{k+1},y,y+z)\frac{(t-s_{k+1})^{n-k}}{(n-k)!}\wt{\nu}^{*k}(dy)\nu(dz)\\
&=&  \sum_{n=0}^\infty
\frac{|\nu|^{n}}{n!}e^{-|\nu|(t-s)}
\int_\Rd \int_\Rd 
\int_s^t  F(v,y,y+z)\sum_{k=0}^n
\frac{n!(v-s)^{k}\wt{\nu}^{*k}(dy)(t-v)^{n-k}}{k!(n-k)!}\nu(dz)dv\\
&=&  \sum_{n=0}^\infty
\frac{|\nu|^{n}}{n!}e^{-|\nu|(t-s)}
\int_\Rd \int_\Rd 
\int_s^t  F(v,y,y+z)
\left((v-s)\wt{\nu}+(t-v)\delta_0\right)^{*n}(dy)\nu(dz)dv\\
&=&  
\int_s^t \int_\Rd \int_\Rd  
F(v,y,y+z) e^{-|\nu|(t-s)}
e^{*\left((v-s)\wt{\nu}+(t-v)\delta_0\right)|\nu|} (dy)\nu(dz)dv\\
&=&  
\int_s^t \int_\Rd \int_\Rd  
F(v,y,y+z) e^{-|\nu|(v-s)}
e^{*(v-s)\nu}(dy)\nu(dz)dv=LHS, \,
\end{eqnarray*}
where  $\delta_0$ is the Dirac measure at $0$.
\qed

In particular, for $s\leq t$ and bounded measurable $F$ we have
\begin{eqnarray}
  \label{eq:d-c}
&& \ex \sum_{s< S_i\leq t}\left[ F(S_i,X_{s,S_{i}-},
  X_{s,S_i})-F(S_i,X_{s,S_{i}-}, X_{s,S_{i}-})\right]\nonumber \\
&& = \ex \int_s^t \int_\Rd \left[F(v,X_{s,v-}, X_{s,v}+z)-F(v,X_{s,v-},X_{s,v-})\right]\nu(dz)dv\,.
\end{eqnarray}
We will consider the filtration
\begin{displaymath}
  \F_t=\sigma\{X_{s,t}\,;\; s\leq t\}\,,\quad t\in \R\,.
\end{displaymath}
For $t\in \R$ we define 
\begin{equation}
  \label{eq:dp}
p_{t}=e^{*t(\nu-|\nu|\delta_0)}
=\sum_{n=0}^\infty \frac{t^n}{n!}(\nu-|\nu|\delta_0)^{*n}
=e^{-t|\nu|}\sum_{n=0}^\infty
\frac{t^n}{n!}\nu^{*n}\,.
\end{equation}
The series converges in the norm of absolute variation of measures.
Clearly, $p_t$ is symmetric,
\begin{equation}
  \label{eq:dpt}
  \frac{\partial}{\partial t}p_t=(\nu-|\nu|\delta_0)*p_t\,,\quad t\in \R\,,
\end{equation}
and $p_{t_1}*p_{t_2}=p_{t_1+t_2}$ for $t_1, t_2\in \R$. 
We have $p_t\geq 0$ for $t\geq 0$, see (\ref{eq:dp}). 
In fact, $p_{u-t}$ is the distribution of $X_{t,u}$, as well as of $X_{t,u-}$, whenever $t\leq u$.
Let 
\begin{equation}
  \label{eq:lke}
  \Psi(\xi)=\int_\Rd (e^{i\xi\cdot z}-1)\nu(dz)\,,\quad \xi \in \Rd\,, 
\end{equation}
where  $\xi\!\cdot\! x$ denotes the usual inner product in $\Rd$. By symmetry of
$\nu$, 
$$\Psi(\xi)=\int_\Rd (\cos \xi \!\cdot\! z-1)\nu(dz)= \Psi(-\xi)\leq 0$$ 
is real valued for all $\xi$. It is also bounded and continuous on $\Rd$.
We have 
\begin{equation}
  \label{eq:lkf}
\hat{p}_t(\xi)=\int_\Rd e^{i\xi\cdot x}p_t(dx)
=e^{t\Psi(\xi)}\,, \quad \xi \in \Rd\,.
\end{equation}
This is the L\'evy-Khinchin formula--a direct
consequence of (\ref{eq:dp})--and $\Psi$ is the corresponding L\'evy-Khinchin
exponent.

Let $g\in \Li$. For $x\in \Rd$, $t\leq u$, we define the {\it
  parabolic}\/ extension of $g$ by 
\begin{displaymath}
  P_{t,u} g(x)=\int_\Rd g(x+y)p_{u-t}(dy)=g*p_{u-t}(x)\,.
\end{displaymath}
This equals $\ex g(x+X_{t,u})$. For $s\leq t \leq u$ we define the
{\it parabolic martingale}
\begin{displaymath}
  G_t=G_t(x;s,u;g)=P_{t,u} g(x+X_{s,t})\,.
\end{displaymath}

\begin{lemma} \label{l:mp}
$G_t$ is a bounded $\{\F_t\}$-martingale on $s\leq t\leq u$.
\end{lemma}
\proof
Independence of increments of $X$ 
%which results from the properties of Poisson random measures, 
yields
\begin{displaymath}
\ex\{g(x+X_{s,u})|\F_t\}=
\ex\{g(x+X_{s,t}+X_{t,u})|\F_t\}=
P_{t,u}g(x+X_{s,t})\,.  \qed
\end{displaymath}
Let $\phi$ be complex-valued and symmetric: $\phi(-z)=\phi(z)$, and let
$|\phi|\leq 1$. 

For $x\in \Rd$, $s\leq t \leq u$, and $f\in \CC$, we define
$F_t=F_t(x;s,u;f,\phi)$ as
\begin{eqnarray*}
&&
\sum_{s< S_i\leq t}
\left[P_{S_i,u}
  f(x+X_{s,S_{i}})-P_{S_i,u}f(x+X_{s,S_{i}-})\right]\phi(X_{s,S_i}-X_{s,S_i-})\\
&& - \int_s^t \int_\Rd \left[P_{v,u}f(x+X_{s,v-}+z)-P_{v,u}f(x+X_{s,v-})\right]\phi(z)\nu(dz)dv\,.
\end{eqnarray*}
\begin{lemma}
  \label{l:fpm}
$\ex |F_t|^p<\infty$ for very $p>0$.
\end{lemma}
\proof
Since $P_{v,u}f$ is bounded for $v\leq u$, the continuous (integral)
part in the definition of $F_t$ is bounded.
We also see that the jump part (the sum above) is bounded by a
constant multiple of $N(s,t)$, which in fact yields exponential integrability of $F_t$.
\qed

In what follows we will denote $\Delta X_{s,t}=X_{s,t}-X_{s,t-}$.
\begin{lemma}
  \label{l:mm}
$\{F_t\}$ is an $\{\F_t\}$-martingale for $s\leq t \leq u$.
\end{lemma}
\proof
By independence of arrivals of signals $\{S_i\}$ on disjoint time
intervals, and by Lemma~\ref{lem:ls},
for $s\leq t_1\leq t_2\leq u$ we have
\begin{eqnarray*}
&&
\ex \left\{\left[ \sum_{t_1< S_i\leq t_2}
\big(P_{S_i,u}f(x+X_{s,S_{i}})-P_{S_i,u}f(x+X_{s,S_{i}-})\big)
\phi(\Delta X_{s,S_i})\right]|\F_{t_1}\right\}\\
&=&
\ex  \sum_{t_1< S_i\leq t_2}
\left[P_{S_i,u}f(x'+X_{t_1,S_{i}})-P_{S_i,u}f(x'+X_{t_1,S_{i}-})\right]
\phi(\Delta X_{t_1,S_i})\\
&=&
\ex \int_{t_1}^{t_2} 
\int_\Rd \left[P_{v,u}f(x'+X_{t_1,v-}+z)-P_{v,u}f(x'+X_{t_1,v-})\right]\phi(z)\nu(dz)dv\,,
\end{eqnarray*}
where $x'=x+X_{s,t_1}$. This gives the martingale property of $F$.
\qed

\begin{lemma}
  \label{l:imf}
$G_t(x;s,u;g)
=F_t(x;s,u;g,1)+P_{s,u}g(x)$.
\end{lemma}
\proof
Since $t\mapsto G_t$ is piecewise differentiable with almost surely
finite number of discontinuities of the first kind (that is,  jumps), we have 
\begin{eqnarray*}
  P_{t,u}g(x+X_{s,t})-P_{s,u}g(x) &=&
\sum_{s<S_i\leq t} 
\left[
P_{S_i,u}g(x+X_{s,S_i})-P_{S_i,u}g(x+X_{s,S_i-})\right]\\
&&+\int_s^t \frac{\partial}{\partial v}P_{v,u}g(x')dv\,,
\end{eqnarray*}
where $x'=x+X_{s,v-}$. This may be considered a version of the It\^o
formula (\cite{MR2020294}). 
The proof is concluded by using (\ref{eq:dpt}),
\begin{eqnarray*}
  \frac{\partial}{\partial v}P_{v,u}g(x')&=&
-\int_\Rd (\nu - |\nu|\delta_0)(dz)P_{v,u}g(x'+z)\\
&=&-\int_\Rd [P_{v,u}g(x+X_{s,v-}+z)-P_{v,u}g(x+X_{s,v-})]\nu(dz)\,. \qed
\end{eqnarray*}
Let $s=t_0\leq t_1 \leq \ldots \leq t_n=t$, and $\sup\{t_i-t_{i-1}\,:\;
i=1,\ldots,n\}\to 0$ as $n\to \infty$.
Since $F_t$ is square integrable, by orthogonality of increments we have  for $s\leq t \leq u$,
\begin{eqnarray*}
\ex F_t^2&=&\ex \sum_{i=1}^n (F_{s,t_i}-F_{s,t_{i-1}})^2\nonumber \\
&\to& 
\ex \sum_{s<S_i\leq t} 
\left[P_{S_i,u}f(x+X_{s,S_i})-P_{S_i,u}f(x+X_{s,S_i-})\right]^2\phi^2(\Delta
X_{s,S_i})\,. \label{eq:qv}
\end{eqnarray*}
The convergence follows from the fact that the integral part of
$F$ is Lipschitz continuous.
Hence the quadratic variation process of $F$ (\cite{MR745449}) is
\begin{equation}
[F,F]_t = \sum_{s<S_i\leq t} 
\left[P_{S_i,u}f(x+X_{s,S_i})-P_{S_i,u}f(x+X_{s,S_i-})\right]^2\phi^2(\Delta
X_{s,S_i})
%\,,\quad s\leq t\leq u
\,. \label{eq:qvp}
\end{equation}
By Lemma~\ref{l:imf}, the quadratic variation of $G$ is 
\begin{equation}
[G,G]_t = |P_{s,u}g(x)|^2+\sum_{s<S_i\leq t} 
\left[P_{S_i,u}g(x+X_{s,S_i})-P_{S_i,u}g(x+X_{s,S_i-})\right]^2
%\,,\quad s\leq t\leq u
\,. \label{eq:qvpg}
\end{equation}
By (\ref{eq:qv}), polarization, and Lemma~\ref{lem:ls}, 
\begin{eqnarray*}
\ex F_t G_t & = & \ex F_t(x;s,u;f,\phi)\left[G_t(x;s,u;g)-P_{s,u}g(x)\right]\\
&=&\ex \sum_{s<S_i\leq t} 
\left[P_{S_i,u}f(x+X_{s,S_i})-P_{S_i,u}f(x+X_{s,S_i-})\right]\\
&&\;\;\;\;\;\;\;\;\;\;\;\;\;\left[P_{S_i,u}g(x+X_{s,S_i})-P_{S_i,u}g(x+X_{s,S_i-})\right]
\phi(\Delta X_{s,S_i})\\
&=&
\ex \int_s^t \int_\Rd 
\left[P_{v,u}f(x+X_{s,v-}+z)-P_{v,u}f(x+X_{s,v-})\right]\\
&&\;\;\;\;\;\;\;\;\;\;\;\;\;\;\;\left[P_{v,u}g(x+X_{s,v-}+z)-P_{v,u}g(x+X_{s,v-})\right]
\phi(z)\nu(dz)dv\\
&=& \int_s^t \int_\Rd p_{v-s}(dy) \int_\Rd
\left[P_{v,u}g(x+y+z)-P_{v,u}g(x+y)\right]\\
&&\;\;\;\;\;\;\;\;\;\;\;\;\;\;\;\;\left[P_{v,u}f(x+y+z)-P_{v,u}f(x+y)\right]\phi(z)\nu(dz)dv\,.
\end{eqnarray*}
By Fubini's Theorem, for any  probability measure $\mu$ and $h\in \Lo$,
\begin{equation}
  \label{eq:wrtlm}
\int_\Rd\int_\Rd h(x+y)\mu(dy)dx=\int h(x)dx\,.
\end{equation}
We define $  |F|_t(x;s,u;f,\phi)$ as
\begin{eqnarray*}
&&
\sum_{s< S_i\leq t} \big[P_{S_i,u}|f|(x+X_{s,S_{i}})+P_{S_i,u}|f|(x+X_{s,S_{i}-})\big]
|\phi|(\Delta X_{s,S_i})\\
&&+
\int_{s}^{t} 
\int_\Rd \big[P_{v,u}|f|(x+X_{s,v-}+z)+P_{v,u}|f|(x+X_{s,v-})\big]|\phi(z)|\nu(dz)dv\,.
\end{eqnarray*}
By Lemma~\ref{lem:ls},
$$
\ex |F|_t= 2 \ex \int_{s}^{t} 
\int_\Rd \big[P_{v,u}|f|(x+X_{s,v-}+z)+P_{v,u}|f|(x+X_{s,v-})\big]\nu(dz)dv\,.
$$
Using (\ref{eq:wrtlm}) we obtain
\begin{eqnarray}
&&\int_\Rd \ex |F|_t(x;s,u;f,\phi)dx \nonumber\\
&=&
 2 \int_\Rd \int_\Rd \nonumber
\int_{s}^{t} \int_\Rd 
\big[P_{v,u}|f|(x+y+z)+P_{v,u}|f|(x+y)\big]\nu(dz)\,dv\, p_{s,v}(dy)\, dx \\  
&=&4 \int_{s}^{t}dv \int_\Rd \nu(dz) \int_\Rd
|f(x)|dx=4(t-s)|\nu|\pl f \pl_1<\infty\label{eq:l1b}
\end{eqnarray}
(compare to Lemma~\ref{l:fpm}). Thus, the following integral is
absolutely convergent
\begin{displaymath}
  I_\phi(f,g)=\int_\Rd \ex F_t(x;s,u;f,\phi) G_t(x;s,u;g) dx\,.
\end{displaymath}
We will now consider $g=f$.
By (\ref{eq:qv}), (\ref{eq:qvpg}) and Lemma~\ref{l:imf}, $F_t(x;s,u;f,\phi)$ is
{\it differentially subordinate}\/ to $G_t(x;s,u;f)$ in that 
$$
0\leq [G,G]_t-[F,F]_t \quad \mbox{\it is non-decreasing for } t\in [s,u]\,.
$$
Therefore, by \cite[Theorem~1]{MR1334160}, we have that
\begin{equation}\label{eq:wmt}
   \ex |F_t(x;s,u;f,\phi)| ^p\leq (p^*-1)^p   \,\ex |G_t(x;s,u;f)|^p\,,\quad s\leq t \leq u\,.
\end{equation}
Here and below we assume (\ref{eq:defp}), in particular $1<p<\infty$.

We note that $G_u(x;s,u;f)=f(x+X_{s,u})$.
Using (\ref{eq:wmt}) and (\ref{eq:wrtlm}) we obtain
\begin{equation} 
  \label{eq:psm}
\int_\Rd \ex  
|F_u(x;s,u;f,\phi)|^p dx 
\leq
(p^*-1)^p   
\int_\Rd \ex |f(x+X_{s,u})|^p dx 
=(p^*-1)^p   
\pl f \pl_p^p\,.
\end{equation}
We consider the linear
functional
$$
\Lq\ni g\mapsto \int_\Rd \ex F_u(x;s,u;f,\phi) g(x+X_{s,u})dx \,.
$$
By H\"older's  inequality, (\ref{eq:psm}) and (\ref{eq:wrtlm}) we have
\begin{equation}
\label{eq:bof}
\int_\Rd \ex |F_u(x;s,u;f,\phi) g(x+X_{s,u})| dx \leq 
(p^*-1)\pl f \pl_p
\pl g \pl_q\,.
%\left[\int_\Rd |f(x)|^p dx\right]^{1/p} \left[\int_\Rd |g(x)|^q dx\right]^{1/q}\,.
\end{equation}
Therefore there is a function $h\in \Lp$ such that 
%KB 8/4 for all $g\in \Lp$
\begin{equation}\label{eq:rf}
  \int_\Rd \ex F_u(x;s,u;f,\phi) g(x+X_{s,u})dx  = \int_\Rd
  h(x)g(x)dx\,,\quad g\in \Lq\,, 
\end{equation}
and 
\begin{equation}
  \label{eq:bm}
\pl h\pl_p\leq (p^*-1)\pl f\pl_p\,.
\end{equation}
We also have that $h\in \Lo$, but
the estimate of $\pl h\pl_1$ depends on $|\nu|$ by (\ref{eq:l1b}).

Consider $\xi\in \Rd$, $e_\xi(x)=e^{i \xi\cdot x}$, and 
%KB 8/4
$\EE_t(x;s,u;\xi)=G_t(x;s,u;e_\xi)$.
To bring about the properties of this martingale we note that by (\ref{eq:lkf})
$$
P_{v,u}e_\xi(x)=\int_\Rd e^{i
  \xi\cdot (x+y)}p_{u-v}(dy)=e_\xi(x)e^{(u-v)\Psi(\xi)}\,, \quad v\leq u \,.
$$
We thus have
\begin{eqnarray*}
  \ex F_t \EE_t &=&
 \int_s^t \int_\Rd p_{v-s}(dy) \int_\Rd
\left[P_{v,u}f(x+y+z)-P_{v,u}f(x+y)\right]e^{(u-v)\Psi(\xi)}\\
&&e^{i\xi\cdot (x+y)}[e^{i \xi\cdot z}-1]\phi(z)\nu(dz)dv\,,
\end{eqnarray*}
hence
\begin{eqnarray*}
I=I_\phi(f,e_\xi)&=&\int_\Rd \int_s^t \int_\Rd \int_\Rd 
\left[P_{v,u}f(x+y+z)-P_{v,u}f(x+y)\right]\\
&&
e^{(u-v)\Psi(\xi)} e^{i\xi\cdot (x+y)}[e^{i \xi\cdot z}-1]\phi(z)\nu(dz) p_{v-s}(dy)dvdx\,.
\end{eqnarray*}
The integral is absolutely convergent by (\ref{eq:l1b}). 
Using (\ref{eq:wrtlm}) and properties of the Fourier transform we obtain 
\begin{eqnarray*}
I&=&
\int_s^t
\int_\Rd 
\int_\Rd \left[P_{v,u}f(x+z)-P_{v,u}f(x)\right]
e^{i\xi \cdot x}dx\,
e^{(u-v)\Psi(\xi)} 
[e^{i \xi \cdot z}-1]\phi(z)\nu(dz)dv\\
&=& 
\int_s^t
\int_\Rd 
\int_\Rd [p_{u-v}*f(x+z)-p_{u-v}*f(x)]e^{i\xi\cdot x}dx\,
e^{(u-v)\Psi(\xi)} [e^{i \xi \cdot z}-1]\phi(z)\nu(dz)dv\\
&=& 
\int_s^t 
\int_\Rd 
[e^{-i\xi\cdot z} 
e^{(u-v)\Psi(\xi)}
\hat{f}(\xi)-e^{(u-v)\Psi(\xi)}\hat{f}(\xi)]
e^{(u-v)\Psi(\xi)} [e^{i \xi \cdot z}-1]\phi(z)\nu(dz)dv\\
&=& 
\hat{f}(\xi)
\int_s^t e^{2(u-v)\Psi(\xi)} dv
\int_\Rd 
|e^{i \xi \cdot z}-1|^2\phi(z)\nu(dz)\,.
\end{eqnarray*}
We have $|e^{i \xi \cdot z}-1|^2=(\cos \xi \!\cdot\! z-1)^2+\sin^2 \xi \!\cdot\! z=2(1-\cos \xi
z)=2 \Re (1-e^{i \xi \cdot z})$. By symmetry of $\phi \nu$,
\begin{displaymath}
  I=
\hat{f}(\xi)\left[e^{2(u-t)\Psi(\xi)}-e^{2(u-s)\Psi(\xi)}\right]\frac{{-}1}{\Psi(\xi)}\int_\Rd
(1-e^{i\xi \cdot z})\phi(z)\nu(dz)\,,\quad \mbox{if } \Psi(\xi)<0\,, 
\end{displaymath}
and $I=0$ if $\Psi(\xi)=0$.
We let $t=u=0$, thus obtaining
$I=\hat{f}(\xi)m_s(\xi)$, where $s<0$, and
\begin{equation}
  \label{eq:fsss}
  m_s(\xi)=\left[1-e^{2|s|\Psi(\xi)}\right]
\frac{\int_\Rd (e^{i\xi \cdot z}-1)\phi(z)\nu(dz)}{\Psi(\xi)}\,,\quad 
\mbox{if } \Psi(\xi)\neq 0\,,
\end{equation}
and $m_s(\xi)=0$ if $\Psi(\xi)=0$.
From (\ref{eq:rf}) applied to $g=e_\xi$ we obtain
\begin{equation}
  \label{eq:m}
  \hat{h}(\xi)=m_s(\xi)\hat{f}(\xi)\,,\quad \xi \in \Rd\,.
\end{equation}
Consider the Fourier multiplier $\M_s$ on $\Lt$ with symbol $m_s$
(bounded by $1$). By (\ref{eq:bm}) the operator uniquely extends to $\Lp$ with norm at most $p^*-1$.
Let 
\begin{eqnarray}
m(\xi)&=&
\frac{\int(e^{i\xi\cdot z}-1)\phi(z)\nu(dz)}{\Psi(\xi)}\nonumber\\
&=&
\frac{\int (\cos \xi \!\cdot\! z-1)\phi(z)\nu(dz)}
{\int (\cos \xi\!\cdot\! z-1)\nu(dz)}
\,,\quad \mbox{ if } 
\Psi(\xi)\neq 0\,   \label{eq:fss}
%
%\Psi(\xi)\neq 0\,,
\end{eqnarray}
and $m(\xi)=0$ if $\Psi(\xi)=\int (\cos \xi \!\cdot\! z-1)\nu(dz)=0$. Clearly, $m=\lim_{s\to -\infty} m_s$, pointwise. 

Let $\M$ be the multiplier on $\Lt$ with symbol $m$.
If $f\in \Lt$, then $\M_s f\to \M f$ in $\Lt$ by Plancherel's Theorem as $s\to-\infty$. 
By Fatou's Lemma and (\ref{eq:bm}) it follows that $\pl \M f\pl_p\leq (p^*-1)\pl f\pl_p$. 
Therefore $\M$ extends uniquely from $\CC$ to $\Lp$ without increasing
the norm, which proves Theorem~\ref{th:bfm} when the L\'evy measure is finite. 

In the general case let $\varepsilon>0$,
and $\nu(B)= V(B\cap \{|x|>\varepsilon\})$.
For every $\xi \in \Rd$, we have that $\cos\xi\!\cdot\! z-1\approx-|z|^2/2$ if
$|z|$ is small. Using (\ref{eq:clm}) we conclude that $m(\xi)$ of
(\ref{eq:fss}) tends to $M(\xi)$ of (\ref{eq:fs}) as $\varepsilon\to 0$. 
The latter is {\it defined}\/ to be zero when its denominator vanishes (see
below in this connection). To complete the proof we use the argument
as in the preceding paragraph 
\qed

We like to remark that an {\it antisymmetric} $\phi$,
$\phi(-z)=-\phi(z)$,  yields zero Fourier symbol in Theorem~\ref{th:bfm}
thus our assumption of symmetry of $\phi$ results in no loss
of generality therein. The case of nonsymmetric $V$, vector-valued
$\phi$, and space-inhomogeneous $V$ and $\phi$ require a further
development of the method presented in this paper.  
%%%%%%%%%%%%%%%%%%%%%%%%%%%
\section{Miscellanea}\label{s:m}
If $\Psi(\xi)=\int_\Rd (\cos \xi\!\cdot\! z -1)V(dz)=0$ for $\xi\neq 0$, then $\supp V\subset A_\xi$, where 
$$A_\xi=\{z\,:\; \xi\!\cdot\! z= 2k\pi \mbox{ for some integer }
k\}\,.$$ 
In particular, $A_\xi$ is discrete in the direction of $\xi$. 
By Fubini's theorem $\{\xi\,:\; \Psi(\xi)=0\}$ has zero Lebesgue
measure. Thus our convention that $M(\xi)=0$ when $\Psi(\xi)=0$,
does not influence the definition of $\M$ on $\Lt$ or $\Lp$.
In fact, $M$ does not generally have a limit where $\Psi(\xi)=0$--the
behavior of (\ref{eq:pm}) at the origin is rather representative here. 
Indeed, assume for simplicity of the discussion that $V$ is finite, compactly
 supported and {\it nondegenerate}, that is not concentrated on a proper
 subspace of $\Rd$. Let $\xi\neq 0$ and assume that $\Psi(\xi)=0$.
 The gradient of $\Psi(\xi)=\int(e^{i\xi\cdot z}-1) \phi(z)V(dz)$ is
 $$  
 i\int_\Rd z e^{i\xi\cdot z} \phi(z)V(dz)=
  i\int_{A_\xi} z \phi(z)V(dz)=0\,,
 $$
 and the Jacobian matrix is $-\int z^T z\, \phi(z)V(dz)$.
 Here $z^T$ denotes the transpose of $z$.
 Thus the first nonzero term in the Taylor expansion of $\Psi(\xi+h)$ at $\xi$ is
 $-\frac{1}{2} \int (z\!\cdot\! h)^2V(dz)<0$ if $h\neq 0$.
 We consider
$$
 \frac
 {-\int_\Rd (z\!\cdot\! h)^2\phi(z)V(dz)}
 {-\int_\Rd (z\!\cdot\! h)^2V(dz)}
 \,.
$$
 The limit of this expression exists if $h=r\eta$, $\eta\in \Rd\setminus\{0\}$, and $r\to 0^+$,
 but in general the limit depends on the direction of $\eta$, compare (\ref{eq:pm}). 
% Thus the definition of {\it continuous} $M$ is
% not always possible.

 \begin{example}\label{ex:hm}
{\rm
We now examine (\ref{eq:pm}). Let $\alpha\in (0,2)$, $j\in
\{1,\ldots,d\}$, and 
$$\mu=\delta_{(1,0,\ldots,0)}+\delta_{(-1,0,\ldots,0)}+
\cdots+\delta_{(0,0,\ldots,1)}+\delta_{(0,0,\ldots,-1)}\,.$$
In polar coordinates we define the L\'evy measure $V(drd\theta)=r^{-1-\alpha}dr
\mu(d\theta)$ (of the symmetric $\alpha$-stable L\'evy
process with independent coordinates \cite{MR1739520}). We have 
\begin{eqnarray}
\Psi(\xi)&=&c_\alpha \int |\xi\!\cdot\!  z|^\alpha \mu(dz) \label{eq:fps}\\
&=&c_\alpha \left( |\xi_1|^\alpha+\cdots+|\xi_d|^\alpha \right)\,,
\nonumber  
\end{eqnarray}
where $c_\alpha=-\pi/(2\sin\frac{\pi\alpha}{2}\Gamma(1+\alpha))$, see
\cite[Chapter 14]{MR1739520}.
Let $\phi(z_1,\ldots,z_d)=1$ if $z_k=0$ for $k\neq j$ and 
%K 8/5
$z_j\neq 0$, 
and let $\phi=0$ otherwise (we observe only the jumps of the first coordinate
process). 
The symbol (\ref{eq:fs}) becomes (\ref{eq:pm}) with $j=1$. 
By Theorem~\ref{th:bfm} the corresponding Fourier
multiplier has norm bounded by $p^*-1$.
%The case of $2$ in the exponents in (\ref{eq:pm}) is obtained 
Letting $\alpha\to 2$ we obtain (\ref{eq:nrt}) by Fatou's Lemma 
(see the end of the proof of Theorem~\ref{th:bfm}). Considering
$\phi=a_j$ on the $j$-th coordinate axis (except
at the origin) for $j=1,\ldots,d$, we conclude that 
\begin{equation} \label{eq:c2rt}
\|\sum_{j=1}^d a_jR_j^2f\|_p\leq (p^*-1)\|f\|_p\,,
\end{equation}
 is valid whenever $|a_j|\leq 1$. 
By considering $\mu$ concentrated on $\sqrt{2}/2(\pm 1,\pm 1)\in \Rt$
and suitably chosen $\phi=\pm 1$ we similarly obtain 
\begin{equation}
  \label{eq:i2rt}
\|2R_jR_kf\|_p\leq (p^*-1)\|f\|_p\,, \quad j\not = k\,.
\end{equation}
 in dimension $d=2$. 
From this, the upper bound $2(p^*-1)$
for the Beurling-Ahlfors operator follows, see Introduction.
}   
\end{example}
\begin{example}
{\rm  Let $d=2$ and $j=1$ in (\ref{eq:pm}). We have
$$
 \left|\frac{\partial}{\partial \xi_1}M(\xi)\right|^2= 
\alpha^2\left[\frac{|\xi_2|^\alpha}{\left(|\xi_1|^\alpha+|\xi_2|^\alpha\right)^2}\right]^2|\xi_1|^{2(\alpha-1)}\,.  
$$
This function is not locally integrable at $\xi_1=0$ if
$0<\alpha<1/2$. Thus the symbol does not satisfy the H\"ormander
condition (\cite{stein}).
}
\end{example}

Denote $M(\xi)=\Psi_\phi(\xi)/\Psi(\xi)$, in (\ref{eq:fs}).
There is a tempered distribution, say $K$, with Fourier transform
$M$, such that $\M \phi=K*\phi$ for smooth compactly supported $\phi$.
It is of interest to represent $\M$ as a limit of integrals.
Let $0<\varepsilon<T<\infty$. 
We will approximate $M$ by
\begin{eqnarray}
&&M_\varepsilon^T(\xi)=\left[e^{\varepsilon\Psi(\xi)}-e^{T\Psi(\xi)}\right]\frac{\Psi_\phi(\xi)}{\Psi(\xi)}
=\int_{\varepsilon}^T \Psi_\phi(\xi) e^{t\Psi(\xi)}dt\nonumber\\
&=&\int_\varepsilon^{T} \Psi_\phi(\xi) e^{t\Psi_\phi(\xi)}e^{t\Psi_{1-\phi}(\xi)}dt
=\int_\varepsilon^{T} \left[\frac{d}{dt} e^{t\Psi_\phi(\xi)}\right]e^{t\Psi_{1-\phi}(\xi)}dt\,,\label{eq:aapm}
\end{eqnarray}
where $\varepsilon\to 0$ and $T\to \infty$ (compare the proof of Theorem~\ref{th:bfm}).
Let $K_\varepsilon^T$ be the (tempered) distribution with Fourier transform $M_\varepsilon^T$.
%\begin{lemma}\label{lem:apm}
%There is a (finite) complex measure $K_s$ on $\Rd$ such that $\widehat{K_s}(\xi)=m_s(\xi)$.   
%\end{lemma}
%\proof
%It follows from (\ref{eq:l1b}) and the discussion following
%(\ref{eq:rf}) that $|K_s|\leq 4|s||\nu|$.
%\qed
If $0\leq \phi\leq 1$, we consider convolution semigroups
$p^\phi_t$ and $p^{1-\phi}_t$ of L\'evy processes with Levy measures
$\phi V$ and $(1-\phi)V$, correspondingly. Motivated by
(\ref{eq:aapm}) we consider
\begin{equation}
  \label{eq:km}
K_\varepsilon^T=\int_\varepsilon^T \left[\frac{d}{dt} p^\phi_t\right]*p_t^{1-\phi}dt\,.
\end{equation}
If $dp^\phi_t/dt$ is a finite measure for $t=\varepsilon$ then it is a
finite measure for all $t\geq \varepsilon$ 
%k 8/7
because $|dp^\phi_t/dt|$ is non-increasing in $t$.
Thus, $K_\varepsilon^T$ is a finite measure and 
$$
K=\lim_{\varepsilon\to 0, T\to \infty}K_\varepsilon^T\,,
$$
as distributions.
In passing we like to note that (\ref{eq:km}) gives an analytic
interpretation to our proof of Theorem~\ref{th:bfm}. 
\begin{example}\label{ex:kernel}
{\rm
When $d=2$, $\alpha=1$ and $j=1$ in (\ref{eq:pm}), the
corresponding multiplier is a singular integral
\begin{equation}
  \label{eq:exsi}
\M f(z)=p.v.\int_\Rt K(z-w)f(w)dw\,,\quad z\in \Rt\,,  
\end{equation}
understood as above, with the kernel
\begin{equation}
  \label{eq:ekm}
  K(x,y)=\frac{-{x^2}+{y^2}+{x^2} \log \big|\frac{x}{y}\big|-{y^2} \log \big|\frac{y}{x}\big|}
              {{{\pi}^2} {{({x^2}-{y^2})}^2}}\,,\quad (x,y)\in \Rt\,.
\end{equation}
To obtain (\ref{eq:ekm}), we denote 
\begin{equation}
  \label{eq:Cd}
  p_t(x)=\frac{1}{\pi}\frac{t}{t^2+x^2}\,,\quad t>0\,,\; x\in \R\,.
\end{equation}
It is the density function of the one-dimensional symmetric
$1$-stable L\'evy (Cauchy) process on the line. 
We have $\widehat{p_t}(\xi)=e^{-t|\xi|}$ for $\xi\in \R$, and 
\begin{equation}
  \label{eq:dpt2}
  \frac{d}{dt}p_t(x)=\frac{1}{\pi}\frac{-t^2+x^2}{(t^2+x^2)^2}\,,
\end{equation}
which is integrable for every $t>0$.
Note that $p_t(x)p_t(y)$, for $(x,y)\in \Rt$, is the transition density of
the Cauchy process with independent coordinates on the plane, compare
Example~\ref{ex:hm}. Our discussion above, (\ref{eq:Cd}) 
and (\ref{eq:dpt2}) yield
$$
K(x,y)=\int_0^\infty
\frac{t (-{t^2}+{x^2})}{{{({t^2}+{x^2})}^2} ({t^2}+{y^2})}dt\,.
$$
Of course, $K(x,y)=K(|x|,|y|)$.
By a change of variable,
% we see that $K$ is homogeneous of degree $-2$: 
\begin{equation}
  \label{eq:hk}
K(hx,hy)=h^{-2}K(x,y)\quad \mbox{if } \quad h>0\,.   
\end{equation}
We will determine $K(1,y)$, where $y>1$.
To this end we observe that
$$
\frac{t (-{t^2}+1)}{{{({t^2}+1)}^2} ({t^2}+{y^2})}
=
\frac{2 t}{{{({t^2}+1)}^2} (-1+{y^2})}
-\frac{t (1+{y^2})}{({t^2}+1) {{(-1+{y^2})}^2}}
+\frac{t (1+{y^2})}{{{(-1+{y^2})}^2}({t^2}+{y^2})}\,.
$$
Integration yields
$$ 
K(1,y)=\frac{-1+{y^2}-(1+{y^2}) \log y}{{{\pi }^2}{{(-1+{y^2})}^2}}\,,
$$ 
and (\ref{eq:ekm}) follows by (\ref{eq:hk}). 
}  
\end{example}
We note a mild singularity of the kernel $K(x,y)$ at $y=0$ in the
previous example, in addition to the usual (critical) singularity at
$(0,0)$ (\cite{stein}). We remark that a stronger singularity may 
be obtained in higher dimensions within the same setup. 
The resulting singularities seem amenable by the Calder\'on-Zygmund theory 
(\cite{MR0100768}), where $L\log L$ integrability and cancellation of
the kernel on the unit sphere are only required to prove the
boundedness of $\M$ on $\Lp$, $1<p<\infty$. 
The emphasis in our paper is, however,  on obtaining good estimates of the
norm of the operator. Also, (\ref{eq:fs}) goes much beyond homogeneous 
symbols (\cite{stein}) and gives a wide and natural class of symbols 
and singular integrals which deserve a further study.
We finally note that the $L^p$ boundedness of
our multipliers may have applications to embedding results for anisotropic
Sobolev spaces as in \cite[Section 2.3]{MR1840499}, \cite[Section 3.1]{JStlp}.

{\bf Acknowledgments.} We thank B.~Dyda and A.~Bielaszewski for discussion and
remarks on the paper. We are grateful to E.M.~Stein for pointing the connection to the Marcinkiewicz multiplier theorem.  We thank two anonymous referees 
for useful comments. The second named author 
gratefully acknowledges 
the hospitality of the Department of Statistics at Purdue
University, where the paper was written in part.

\bibliographystyle{abbrv}
\bibliography{martingale_transforms}

%                            ---------- o ----------

\end{document}